\documentclass[smallextended,envcountsect]{svjour3} 
\smartqed 
\usepackage{graphicx}
\usepackage{color, bm, amsmath, amsfonts, url, amssymb, verbatim} 
\usepackage{longtable, multirow}
\usepackage{float}
\usepackage{tikz}
\usepackage[ruled, lined]{algorithm2e}
\newcommand{\norm}[1]{||#1||}
\newcommand{\projP}[2]{\mathcal{P}(#1 \, ; \, #2)}
\newcommand{\projPi}[3]{p_{#1}(#2 \, ; \, #3)}
\newcommand{\mbf}[1]{\mathbf{#1}}
\DeclareSymbolFont{largesymbolsA}{U}{txexa}{m}{n}
\DeclareMathSymbol{\varprod}{\mathop}{largesymbolsA}{16}
\journalname{JOTA}

\begin{document}

\title{A fast and convergent combined Newton and gradient descent method for computing steady states of chemical reaction networks}


\author{Silvia Berra \and  Alessandro La Torraca \and Federico Benvenuto \and Sara Sommariva}

\institute{Silvia Berra \at
           Dipartimento di Matematica, Universit\`a di Genova, via Dodecaneso 35 16146 Genova, Italy \\
           silvia.berra@dima.unige.it  
           \and
           Alessandro La Torraca \at
           Data and Analytics Chapter, Roche S.p.A., Monza, Italy\\
           alessandro.la\_torraca@roche.com 
           \and
           Federico Benvenuto  \at
           Dipartimento di Matematica, Universit\`a di Genova, via Dodecaneso 35 16146 Genova, Italy \\
           benvenuto@dima.unige.it
            \and
           Sara Sommariva,  Corresponding author  \at
           Dipartimento di Matematica, Universit\`a di Genova, via Dodecaneso 35 16146 Genova, Italy \\
           sommariva@dima.unige.it 
}


\maketitle

\begin{abstract}
In this work we present a fast, globally convergent, iterative algorithm for computing the asymptotically stable states of nonlinear large--scale systems of quadratic autonomous Ordinary Differential Equations (ODEs) modeling, e.g., the dynamic of complex chemical reaction networks. 
Towards this aim, we reformulate the problem as a box--constrained optimization problem where the roots of a set of nonlinear equations need to be determined. 
Then, we propose to use a projected Newton's approach combined with a gradient descent algorithm so that every limit point of the sequence generated by the overall algorithm is a stationary point. 
More importantly, we suggest replacing the standard orthogonal projector with a novel operator that ensures the final solution to satisfy the box constraints while lowering the probability that the intermediate points reached at each iteration belong to the boundary of the box where the Jacobian of the objective function may be singular.
The effectiveness of the proposed approach is shown in a practical scenario concerning a chemical reaction network modeling the signaling network of colorectal cancer cells. 
Specifically, in this scenario the proposed algorithm is proven to be faster and more accurate than a classical dynamical approach where the asymptotically stable states are computed as the limit points of the flux of the Cauchy problem associated with the ODEs system.
\end{abstract}

\section{Introduction}
\label{par2}

This paper is focused on the solution of root-finding problems in several variables where the system is composed by algebraic second degree equations. This kind of problems are of interest in many application areas, including queuing problems, neutron transport theory, linear quadratic differential games (see e.g. \cite{poloni2013quadratic} and references therein).
Our work is motivated by the study of non--negative steady states of biological interaction networks which frequently arise in systems biology \cite{feinberg1995,conradi2005,gabor2015}.
In particular, the problem of determining the steady states of complex Chemical Reaction Networks (CRNs) in healthy and cancer cells is considered \cite{Jordan,Sever}. 
Indeed, by applying the law of mass action, the kinetics of the concentration of the proteins involved in the network can be modelled by a large first order polynomial system of Ordinary Differential Equations (ODEs).
Finally, when no exogenous factors are considered, the equation system is quadratic and autonomous \cite{Feinberg,Yu_Craciun_2018,Chellaboina}.
From an abstract point of view, steady states of large system of quadratic equation are far from being known, as a general theory exists up to the two dimensional case \cite{Reyn2007PhasePO}.
In the case of an ODE representing a CRN, unknown concentrations cannot assume negative values, and then the asymptotically steady states must fulfill the non--negative constrained algebraic system of equations deriving from setting equal to zero the time derivatives of the ODE system.
The number of involved unknown protein concentrations may scale up to several hundreds and an efficient and accurate algorithm to solve the non--negative steady state problem is at the basis of tuning the kinetic parameters of the ODE system starting from experimental data, thus enabling the study of cell cancer behaviour in real applications.

From a computational point of view, equilibrium can be found either in the direct way by taking the limit of the flux of the ODEs, or imposing the vanishing of the derivatives and solving the corresponding root-finding problem.
The direct approach is computationally expensive, especially when the orbits of the dynamical system are bent around the equilibrium point as the time to run across the orbit may become arbitrarily large \cite{DORMAND198019}.
On the other hand, the main drawback of the second strategy is that these systems do not usually show any mathematical property which can ensure convergence of the root-finding algorithm.
Indeed, pertinent good mathematical properties, such as matrix positive definiteness, depend on the form of the considered biological network, and are not ensured in a general case.
The typical structure of the ODE system associated with a chemical reaction network (based on the mass action law) not only prevents us from exploiting recent methods to find the steady state solutions by solving vector quadratic equations \cite{poloni2013quadratic}, but it also makes it difficult to use classical methods, such as the Newton's or gradient descent methods.
Indeed, as the non--negative steady state normally belongs to the frontier of the positive cone and therefore it has many components equal to zero, classical non--negative projected Newton-type methods are unstable as the Jacobian matrix - computed in a neighborhood of the solution - is strongly sparse and non-invertible.
Moreover, classical projected gradient methods are known to be stable but with slow convergence, especially in cases of coupling them with a non--negative projection.

In this work we propose to overcome these limitations, by introducing a root-finding strategy based on combining steps of Newton's method and steps of gradient descent.
While the Newton's method is applied to the algebraic equation system, the gradient descent is applied by scalarizing the system, i.e. minimizing the norm of the l.h.s. of the equation system \cite{Khanh1993OptimalityCV}.
To make the Newton's method more stable, instead of the standard orthogonal projection, we use a non--linear projection operator onto the non--negative orthant that is substantially a idempotent operator providing small positive entries rather than zero components.
Doing so, it improves the condition number of the Jacobian matrix preventing the Newton's step to be unstable and hence making regularization unnecessary.
Therefore we combine the (non--linearly projected) Newton's method with a gradient method to iteratively refine the starting point of the former until we get the convergence to a non--negative stationary point.
We prove the convergence of this combined technique provided that a proper backtracking rule on the gradient method is considered.
Moreover, we test the efficiency of the proposed technique in the case of simulated CRN data, showing that, compared to standard ODE solvers, this method computes the steady states achieving greater accuracy in less time.
The MATLAB\textsuperscript{\textregistered} codes implementing the proposed approach is freely available at the GitHub repository \url{https://​github.​com/​theMI​DAgro​up/​CRC_​CRN.git}. 

The rest of the paper is organized as follows. 
In Sect. \ref{sec:math_pb} we introduce the mathematical formulation of the problem and we describe the proposed algorithm whose converge properties are studied in Sect. \ref{sec:NLPC}. 
In Sect. \ref{sec:NLPC_for_CRN} we consider the problem of finding the asymptotically steady states of a CRN and we reformulate it as a non--negative constrained root--finding problem.
In Sect. \ref{sec:results} we show the results obtained by applying NLPC to a CRN designed for modeling cell signaling in colorectal cells and the most common mutations occurring in colorectal cancer. 
Finally our conclusions are offered in Sect. \ref{sec:conclusions}.

\section{Mathematical formulation}
\label{sec:math_pb}
We consider the box--constrained set of nonlinear equations
\begin{equation}\label{eq:box_eqs}
\begin{cases}
\mbf{f}(\mbf{x}) = \mbf{0} \\
\mbf{x} \in \Omega
\end{cases}
\end{equation} 
where $\Omega = \varprod_{i=1}^n \Omega_i \subseteq \mathbb{R}^n$ is the Cartesian product of $n$ closed intervals $\Omega_i \subseteq \mathbb{R}$, and $\mbf{f}: \mathbb{R}^n \rightarrow \mathbb{R}^n$ is a continuously differentiable function on $\Omega$. In the considered problem of finding the non--negative steady states of quadratic autonomous ODEs systems, $\mbf{f}$ is composed by second--degree polynomials and $\Omega$ is the positive convex cone.

Several numerical approaches have been proposed to solve (\ref{eq:box_eqs}). Among these, a classical fast approach is the projected Newton's method \cite{Bertsekas1982,bertsekas1997} where the projector on the closed convex set $\Omega$,  $P:\mathbb{R}^n \rightarrow \Omega$ such that for all $\mbf{z} \in \mathbb{R}^n$
\begin{equation}\label{eq:classical_P}
P(\mbf{z}) = \underset{\mbf{y} \in \Omega}{\mathrm{argmin}} \norm{ \mbf{y} - \mbf{z}} \, ,
\end{equation}
is applied at each iteration of a Newton's scheme so that the final solution satisfies the box constraints in (\ref{eq:box_eqs}). 
However, in the general case convergence properties of the projected Newton's methods strongly depend on the initial point, as no global convergence is guaranteed \cite{nesterov2006}. 
Additionally, the standard orthogonal projector $P$ tends to provide iterative estimates on the boundary of $\Omega$ (e.g. when $\Omega$ is the positive cone $P$ sets to zero all the negative components) and therefore it may compromise the stability of the Newton's method as the Jacobian of $\mbf{f}$ can be singular computed at these boundary estimates. 

An alternative approach to Newton's method consists in using the projected gradient descent method \cite{goldstein1964,levitin1966} for solving the optimization problem 
\begin{equation}\label{eq:eqs_min}
\mbf{x} = \underset{\mbf{x} \in \Omega}{\mathrm{arg min}} \, \Theta(\mbf{x}) \, ,
\end{equation}
where 
\begin{equation}\label{eq:def_theta}
    \Theta(\mbf{x}) = \frac{1}{2} \norm{\mbf{f}(\mbf{x})}^2 \, .
\end{equation}
As opposite to Newton's method, many convergence results may be proved for the projected gradient methods, see e.g. \cite{bertsekas1997,wang2000} and references therein. On the other hand, the projected gradient method only has a sub--linear convergence rate and thus results to be slower than the Newton's algorithm also when properly designed strategies for selecting the stepsize are used \cite{barzilai1988,crisci2019,dai2005,serafini2005}. 

Motivated by this consideration, some recent works have proposed to combine the two approaches \cite{shi1996,han2003,chen2017,di2021}. Along this lines we present the Non--Linearly Projected Combined (NLPC) method that is summarized in Algorithm \ref{algo:gp_newton}. The main ideas behind NLPC are two. First of all, we replace the classical projector with a novel operator $\mathcal{P}$, introduced in the following definition, that ensures the constraint $\mbf{x} \in \Omega$ to be respected while lowering the probability that the points defined at each iteration reach the boundary of $\Omega$.

\begin{definition}\label{def:our_projector}
Given $\Omega = \varprod_{i=1}^n \Omega_i$, $\Omega_i \subseteq \mathbb{R}$ convex for all $i \in \{1, \dots, n\}$, and given $\mbf{x} = \left(x_1, \dots, x_n \right)^\top \in \Omega$, we define the operator $\projP{\, \cdot}{\mbf{x}}: \mathbb{R}^n \rightarrow \Omega$, so that, for all $\mbf{z} \in \mathbb{R}^n$, $\projP{\mbf{z}}{\mbf{x}} = \left(\projPi{1}{z_1}{x_1}, \dots, \projPi{n}{z_n}{x_n} \right)^\top$ where
\begin{equation*}
\projPi{i}{v}{w} = 
\begin{cases}
v \quad \text{if} \quad v \in \Omega_i\\
w \quad \text{if} \quad v \not\in \Omega_i
\end{cases}
\end{equation*}
\label{defP}
with $v\in\mathbb R$ and $w\in\Omega_i \subset \mathbb R$.
\end{definition}

The second idea behind NLPC method was inspired by \cite{chen2017} and consists in trying at each iteration a fixed number of step lengths $\alpha^j,\, j \in \{ 0, \dots, J\}$, along the Newton's direction $\mbf{d}_k$, where $\mbf{d}_k$ is defined as the solution of the set of equations $\mathbf{J}_{\mbf{f}}(\mbf{x}_k) \mbf{d}_k = -\mbf{f}(\mbf{x}_k)$, being $\mathbf{J}_{\mbf{f}}(\mbf{x}_k)$ the Jacobian matrix of $\mbf{f}$ evaluated in $\mbf{x}_k$. If none of the tested stepsizes satisfies the Armijo Rule 
\begin{equation}\label{eq:AR_newton}
          \|\mbf{f}(\projP{\mbf{x}_k + \alpha^{j} \mbf{d}_k}{\mbf{x}_k} )\| \leq \sqrt{1-  \alpha^{j} \sigma_N} \ \|\mbf{f}(\mbf{x}_k)\| \, ,
\end{equation}
we then move along the gradient descent direction with a stepsize chosen so as to satisfy two conditions that, as we shall prove in Theorem \ref{th:conv_analysis}, guarantee a convergence result for NLPC algorithm.

\begin{algorithm}
\SetKwInOut{Input}{Input}\SetKwInOut{Output}{Output}\SetKwInOut{AND}{and}
\SetKw{And}{\hspace{\algoskipindent}\itshape and\;}
\SetKwBlock{Condition}{}{}
\DontPrintSemicolon 
\Input{$\mbf{x}_0 \in \Omega$; $\tau \in (0, +\infty)$; $\alpha $, $\sigma_N$, $\sigma_G$, $\rho \in (0, 1)$; $J \in \mathbb{N} \setminus \{0\}$} 
$ FLAG \leftarrow 0$; $ k \leftarrow 0$\;
\While{$\|\mbf{f}(\mbf{x}_k)\| > \tau$}{
     \eIf{FLAG=0}{solve $\mathbf{J}_{\mbf{f}}(\mbf{x}_k) \mbf{d}_k = -\mbf{f}(\mbf{x}_k) $\;
      $j \leftarrow 0$\; 
      \While{$j\ \leq\ J$ }{
       $\mbf{x}_{k+1} = \projP{\mbf{x}_k + \alpha^{j} \mbf{d}_k}{\mbf{x}_k}$ \;
      \eIf{$\|\mbf{f}(\mbf{x}_{k+1} )\| \leq \sqrt{1-  \alpha^{j} \sigma_N} \ \|\mbf{f}(\mbf{x}_k)\|$}{
      $j \leftarrow J + 1$; $FLAG \leftarrow 0$; $k \leftarrow k + 1$}{
     $j \leftarrow j+1$;  $ FLAG \leftarrow 1$\;}
      }
      }{
      $\mbf{d}_k = - \nabla \Theta(\mbf{x}_k)$ \;
      $j \leftarrow 0$\; 
      \While{$FLAG = 1$}{
       $\mbf{x}_{k+1} = \projP{\mbf{x}_k + \alpha^j \mbf{d}_k}{\mbf{x}_k} $ \;
      \eIf{\Condition{\mbox{$\Theta(\mbf{x}_{k+1}) \leq \Theta(\mbf{x}_k) + \sigma_G \nabla \Theta(\mbf{x}_k)^T (\mbf{x}_{k+1}  - \mbf{x}_k )$} \;
      \And
      \mbox{$\sqrt{\textstyle\sum\limits_{i \in \mathcal{M}_{\alpha^j}(\mbf{x}_k)} (P_i(x_{k,i} + d_{k,i}) - x_{k,i})^2 } \geqslant \rho \sqrt{\textstyle\sum\limits_{i \in \mathcal{N}_{\alpha^j}(\mbf{x}_k)} (P_i(x_{k,i} + d_{k,i}) - x_{k,i})^2 }$} \;}}{
      $FLAG \leftarrow 0$; $k \leftarrow k + 1$ \;
      }{
      $j \leftarrow j+1$}
      }
      }
    }
\caption{The NLPC algorithm}
\label{algo:gp_newton}
\end{algorithm}

\section{Convergence properties of the NLPC method}\label{sec:NLPC}

Now we present a convergence analysis of NLPC algorithm after describing the main tools exploited in the algorithm and the main properties of $\mathcal{P}$.

\begin{definition}\label{def:set_B_M_N}
Given $\mbf{x} \in \Omega$, $\mbf{d} \in \mathbb{R}^n \setminus \{\mbf{0}\} $ and $\alpha > 0$, we define 
\begin{equation}\label{def:b_xd}
    \mathcal{B}(\mbf{x}, \mbf{d}) := \left\{ i \in \{1, \dots, n\}\ s.t.\ x_i + \alpha d_i  \notin \Omega_i \ \forall \alpha>0 \right\}
\end{equation}
\begin{equation}\label{def:r_xd}
    \mathcal{M}_{\alpha}(\mbf{x}, \mbf{d}) := \left\{ i \in \{1, \dots, n\}\ s.t.\ x_i + \alpha d_i  \in \Omega_i  \right\}
\end{equation}
\begin{equation}\label{def:c_xd}
    \mathcal{N}_{\alpha}(\mbf{x}, \mbf{d}) := \{1, \dots, n\}\setminus \left( \mathcal{B}(\mbf{x}, \mbf{d}) \cup \mathcal{M}_{\alpha}(\mbf{x}, \mbf{d})  \right)
\end{equation}
For the ease of notation, when $\mbf{d} = - \nabla\Theta(\mbf{x})$, the set defined  in (\ref{def:b_xd}), (\ref{def:r_xd}) and (\ref{def:c_xd}) will be simply denoted as $\mathcal{B}(\mbf{x})$, $\mathcal{M}_{\alpha}(\mbf{x})$, and $\mathcal{N}_{\alpha}(\mbf{x})$, respectively.
\end{definition}

\begin{figure}[H]
\begin{center}
\includegraphics[width=12cm]{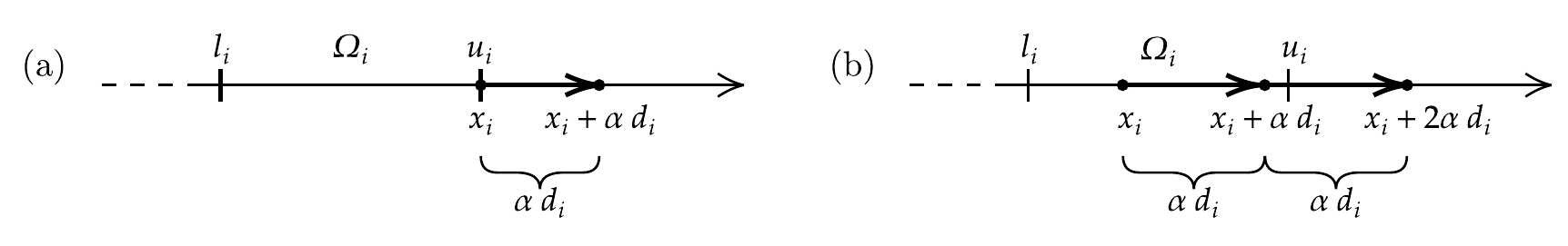}
\end{center}
\caption{(a) Example where $i \in \mathcal{B}(\mbf{x}, \mbf{d})$. (b) Example where $i \in \mathcal{M}_{\alpha}(\mbf{x}, \mbf{d})$ and $i \in \mathcal{N}_{2 \alpha}(\mbf{x}, \mbf{d})$. }\label{fig:sets_1D}
\end{figure}

\begin{remark}
It can be easily shown that, for all $\alpha > 0$,
\begin{equation}
\mathcal{B}(\mbf{x}, \mbf{d})\ \cup \  \mathcal{M}_{\alpha}(\mbf{x}, \mbf{d}) \ \cup \ \mathcal{N}_{\alpha}(\mbf{x}, \mbf{d}) \ = \ \left\{1, \dots, n \right\} ,  
\label{sets_BMN}
\end{equation}
and the three sets are pairwise disjoint. More in detail, as illustratively depicted in Figure \ref{fig:sets_1D}(a), the set $\mathcal{B}(\mbf{x}, \mbf{d})$ contains all the coordinates $i$ that prevent $\mbf{d}$ from being a feasible direction as moving along the corresponding component $d_i$ violates the constraint in (\ref{eq:box_eqs}). As an example, when $\Omega = \varprod_{i=1}^n [\ell_i, u_i]$, $\ell_i < u_i$,
\begin{displaymath}
\mathcal{B}(\mbf{x}, \mbf{d}) = \left\{ i \in \{1, \dots, n\}\ s.t. \ (x_i = \ell_i \land d_i < 0) \lor (x_i = u_i \land d_i > 0)\right\} \ .
\end{displaymath}
\end{remark}
Instead, fixed a stepsize $\alpha>0$, $\mathcal{M}_{\alpha}(\mbf{x}, \mbf{d})$ contains all the components $i$ for which $x_i + \alpha d_i$ still satisfies the constraint of the problem, while $\mathcal{N}_{\alpha}(\mbf{x}, \mbf{d})$ collects the components for which the stepsize $\alpha$ is too big, but a feasible vector may be found by lowering it, see Figure \ref{fig:sets_1D}(b). \\

\begin{proposition}\label{prop:proprieties_P}
Given $\mbf{x} \in \Omega$ and $\mbf{d} \in \mathbb{R}^n$, it holds
\begin{description}
\item[(a)] For all $\alpha > 0$
\begin{equation}\label{eq:prop_prod}
    \left(\mbf{x} - \projP{\mbf{x}+\alpha\mbf{d}}{\mbf{x}} \right)^T \left(\mbf{x} + \alpha  \mbf{d} - \projP{\mbf{x}+\alpha\mbf{d}}{\mbf{x}} \right) = 0
\end{equation}
and 
\begin{equation}\label{eq:prop_norm}
    \norm{\projP{\mbf{x}+\alpha\mbf{d}}{\mbf{x}} - \mbf{x} } = \alpha \sqrt{\sum_{i \in \mathcal{M}_{\alpha}(\mbf{x}, \mbf{d})}{d_i^2}} \ . 
\end{equation}
\item[(b)] $\mathbf{g} : (0, \infty) \rightarrow \Omega$ s.t. $\mathbf{g}(\alpha) = \projP{\mbf{x} + \alpha \mbf{d}}{\mbf{x}}$ is continuous in 0.
\item[(c)] $\varphi : (0, \infty) \rightarrow \mathbb{R}$ s.t. $\varphi(\alpha) = \frac{\norm{\projP{\mbf{x}+\alpha\mbf{d}}{\mbf{x}}- \mbf{x}}}{\alpha}$ is monotonically nonincreasing. \label{phi_decreases}
\end{description}
\end{proposition}
\begin{proof}
(a) Equations (\ref{eq:prop_prod}) and (\ref{eq:prop_norm}) follow from Definition \ref{def:our_projector} which implies 
\begin{displaymath}
\begin{split}
    ( \mbf{x} -  & \projP{\mbf{x} + \alpha \mbf{d}}{\mbf{x}} )^T \left(\mbf{x} + \alpha  \mbf{d} - \projP{\mbf{x} + \alpha \mbf{d}}{\mbf{x}} \right) = \\
    & \sum_{i=1}^n {\left( x_i -  \projPi{i}{x_i + \alpha d_i}{x_i} \right) \left(x_i + \alpha  d_i -  \projPi{i}{x_i + \alpha d_i}{x_i}) \right)} = 0 \, 
\end{split}
\end{displaymath}
and 
\begin{displaymath}
\begin{split}
 \norm{\projP{\mbf{x}+\alpha\mbf{d}}{\mbf{x}} - \mbf{x} } & = \sqrt{\sum_{i=1}^n (\projPi{i}{x_i + \alpha_i d_i}{x_i}-x_i)^2} = \alpha \sqrt{\sum_{i\in\mathcal{M}_{\alpha}(\mbf{x}, \mbf{d})} d_i^2} \, .
 \end{split}
\end{displaymath}
(b) The result directly follows from equation (\ref{eq:prop_norm}). Indeed
\begin{displaymath}
\norm{\mathbf{g}(\alpha) - \mathbf{g}(0)} = \norm{\projP{\mbf{x}+\alpha\mbf{d}}{\mbf{x}} - \mbf{x} } = \alpha \sqrt{\sum_{i \in \mathcal{M}_{\alpha}(\mbf{x}, \mbf{d})}{d_i^2}} \leq  \alpha  \  \norm{\mbf{d}} \xrightarrow[\alpha \to 0^+]{} 0 \, .
\end{displaymath}
(c) We observe that equation (\ref{eq:prop_norm}) implies $\varphi(\alpha) = \sqrt{\sum_{i\in \mathcal{M}_{\alpha}(\mbf{x}, \mbf{d})} d_i^2}$.
Since $\Omega_i$ is a convex set, given $0 \leq \alpha_1 \leq \alpha_2$ it holds $\mathcal{M}_{\alpha_2}(\mbf{x}, \mbf{d}) \subseteq \mathcal{M}_{\alpha_1}(\mbf{x}, \mbf{d})$, and thus
\begin{displaymath}
\begin{split}
\varphi(\alpha_1) - \varphi(\alpha_2) & = \sqrt{\sum_{i\in \mathcal{M}_{\alpha_1}(\mbf{x}, \mbf{d})} d_i^2} - \sqrt{\sum_{i\in \mathcal{M}_{\alpha_2}(\mbf{x}, \mbf{d})} d_i^2} \geq 0
\end{split}
\end{displaymath}
\qed
\end{proof}

As we shall see in the next theorems, the results shown in Proposition \ref{prop:proprieties_P} allow us to prove convergence properties of the proposed NLPC algorithm similar to those holding when the classical projector on the closed set $\Omega$ is employed instead of the operator $\mathcal{P}$  \cite{bertsekas1997,chen2017}.

\begin{theorem}\label{th:stat_point}
Given $\Theta : \mathbb{R}^n \rightarrow \mathbb{R}$ a continuously differentiable function on $\Omega$ and $\mbf{x} \in \Omega$, then $\mbf{x}$ is a stationary point of $\Theta$ in $\Omega$ iff
\begin{equation}\label{eq:proj_and_stat}
\projP{\mbf{x} - \alpha \nabla\Theta(\mbf{x})}{\mbf{x}} = \mbf{x} \quad \forall \, \alpha > 0 \, .
\end{equation}
\label{th:stationarity}
\end{theorem}
\begin{proof} 
Let's consider the projector P on the closed convex set $\Omega$, defined in Eq. (\ref{eq:classical_P}).
 The following properties hold: (i) $P(\mbf{z}) = (P_1(z_1), \dots, P_n(z_n))$, being 
\begin{displaymath}
P_i(z_i) = \left\{
\begin{array}{cl}
    \ell_i & \quad \textrm{if} \quad  z_i \leq \ell_i \\
    z_i & \quad \textrm{if} \quad \ell_i < z_i < u_i \\
    u_i & \quad \textrm{if} \quad z_i \geq u_i
\end{array}\right.  
\end{displaymath}
where we denoted $\Omega_i = [\ell_i, u_i]$ with $\ell_i,\, u_i \in \mathbb{R} \cup \left\{\pm \infty \right\}$; and (ii) $\mbf{x}$ is a stationary point iff $P(\mbf{x} - \alpha \nabla\Theta(\mbf{x})) = \mbf{x}$ $\forall \, \alpha > 0$ \cite{bertsekas1997}.

We now assume that condition (\ref{eq:proj_and_stat}) holds and thus, $\forall i \in \{1, \dots, n\}$,
\begin{displaymath}
\projPi{i}{x_i - \alpha \partial_i\Theta(\mbf{x})}{x_i} = x_i \quad \forall \alpha > 0 \, .
\end{displaymath}
For each $i \in \{1, \dots, n\}$, we then have only three possibilities:
\begin{itemize}
    \item $x_i \in (\ell_i, u_i)$ and $\partial_i\Theta(\mbf{x}) = 0$. Then $P_i(x_i - \alpha \partial_i\Theta(\mbf{x})) = P(x_i) = x_i \ \forall \alpha > 0$;
    \item $x_i = \ell_i$ and $\partial_i\Theta(\mbf{x}) \geq 0$. In this case,
    $P_i(x_i - \alpha \partial_i\Theta(\mbf{x})) = \ell_i = x_i \ \forall \alpha > 0$
    \item $x_i = u_i$ and $\partial_i\Theta(\mbf{x}) \leq 0$. In this case, $P_i(x_i - \alpha \partial_i\Theta(\mbf{x})) = u_i = x_i \ \forall \alpha > 0$
\end{itemize}
In all three cases we obtained $P_i(x_i - \alpha \partial_i\Theta(\mbf{x})) = x_i$ $\forall \alpha > 0$. This implies that $\mbf{x} \in \Omega$ is a stationary point on $\Omega$. 

Conversely, consider a stationary point $\mbf{x} \in \Omega$ and let's assume it exists $\alpha > 0$ such that $\mathcal{P}(\mbf{x} - \alpha \nabla\Theta(\mbf{x}); \mbf{x}) \neq \mbf{x}$. From Proposition \ref{prop:proprieties_P} (a) it follows
\begin{displaymath}
\begin{split}
    0 & =  \left(\mbf{x} - \projP{\mbf{x}-\alpha\nabla\Theta(\mbf{x})}{\mbf{x}} \right)^T \left(\mbf{x} -\alpha\nabla\Theta(\mbf{x}) - \projP{\mbf{x}-\alpha\nabla\Theta(\mbf{x})}{\mbf{x}} \right) \\
    & = \norm{\mbf{x} -  \projP{\mbf{x}-\alpha\nabla\Theta(\mbf{x})}{\mbf{x}}}^2 - \alpha \nabla\Theta(\mbf{x})^T \left(\mbf{x} - \projP{\mbf{x}-\alpha\nabla\Theta(\mbf{x})}{\mbf{x}} \right) 
\end{split}
\end{displaymath}
and thus 
\begin{equation}
  \nabla\Theta(\mbf{x})^T \left( \projP{\mbf{x}-\alpha\nabla\Theta(\mbf{x})}{\mbf{x}} - \mbf{x}\right) = - \frac{\norm{\mbf{x} -  \projP{\mbf{x}-\alpha\nabla\Theta(\mbf{x})}{\mbf{x}}}^2}{\alpha} < 0 \, .
  \label{eq9} 
\end{equation}
Equation (\ref{eq9}) contradicts the assumption of $\mbf{x}$ being a stationary point, that would imply $\nabla\Theta(\mbf{x})^T(\mbf{z}-\mbf{x}) \geq 0$ $\forall \, \mbf{z} \in \Omega$ \cite{bertsekas1997}.\\
\qed
\end{proof}

\begin{theorem}\label{th:descent_dir}
Given $\Theta : \mathbb{R}^n \rightarrow \mathbb{R}$ a continuously differentiable function on $\Omega$  and $\mbf{x} \in \Omega$ that is not a stationary point of $\Theta$, then it exists $\alpha^*>0$ such that $\forall \alpha \in (0,\alpha^*]$ $\left( \projP{\mbf{x} - \alpha \nabla \Theta(\mbf{x})}{\mbf{x}} - \mbf{x} \right)$ is a descent direction for $\Theta$.
\end{theorem}
\begin{proof}
Since $\mbf{x}$ is not a stationary point, according to Theorem \ref{th:stat_point} it exists $\alpha^* > 0$ such that $\mathcal{P}(\mbf{x} - \alpha^* \nabla\Theta(\mbf{x}); \mbf{x}) \neq \mbf{x}$ and thus $\mathcal{P}(\mbf{x} - \alpha \nabla\Theta(\mbf{x}); \mbf{x}) \neq \mbf{x} \ \forall \alpha \in (0, \alpha^*]$ because $\Omega_i$ is a convex set $\forall \ i \in \{1, \dots, n\}$. Therefore, from (\ref{eq9}) it follows $\nabla\Theta(\mbf{x})^T \left( \projP{\mbf{x}-\alpha\nabla\Theta(\mbf{x})}{\mbf{x}} - \mbf{x}\right) < 0$.
\qed
\end{proof}

\begin{theorem}\label{th:GC}
Given $\Theta : \mathbb{R}^n \rightarrow \mathbb{R}$ a continuously differentiable function on $\Omega$, $\mbf{x} \in \Omega$ and $\sigma_G \in (0, 1)$, it exists $\overline{\alpha}>0$ so that for all $\alpha \in (0, \overline{\alpha}]$
\begin{equation}\label{AR_gradient}
    \Theta(\projP{\mbf{x}-\alpha \nabla \Theta(\mbf{x})}{\mbf{x}}) \leqslant \Theta(\mbf{x}) + \sigma_G \nabla \Theta(\mbf{x})^T \left(\projP{\mbf{x} - \alpha \nabla \Theta(\mbf{x})}{\mbf{x}} - \mbf{x} \right) \, .
\end{equation}
\end{theorem}

\begin{proof}
If $\projP{\mbf{x}-\alpha \nabla \Theta(\mbf{x})}{\mbf{x}} = \mbf{x}$ for all $\alpha>0$, then the thesis holds for any $\overline{\alpha} > 0$. Therefore we can assume it exists $\widetilde{\alpha} \in (0, 1)$ such that $\projP{\mbf{x}-\alpha \nabla \Theta(\mbf{x})}{\mbf{x}} \neq \mbf{x}$ for all $\alpha \in (0, \widetilde{\alpha}]$. In the following we shall denote $\mbf{x}(\alpha) := \projP{\mbf{x}-\alpha \nabla \Theta(\mbf{x})}{\mbf{x}}$. \\ By the mean value theorem, it exists $\boldsymbol{\xi}_{\alpha}$ on the segment between $\mbf{x}$ and $\mbf{x}(\alpha)$ so that
\begin{displaymath}
\begin{split}
    \Theta(\mbf{x}(\alpha))  - \Theta(\mbf{x}) & = 
     \nabla\Theta(\boldsymbol{\xi}_{\alpha})^T \left(\mbf{x}(\alpha) - \mbf{x}\right) \\
    & = \sigma_G \nabla\Theta(\mbf{x})^T \left( \mbf{x}(\alpha) -\mbf{x} \right)  - (\sigma_G -1) \nabla\Theta(\mbf{x})^T \left(\mbf{x}(\alpha) - \mbf{x} \right) + \\
    & \qquad + \left(\nabla\Theta(\boldsymbol{\xi}_{\alpha}) - \nabla\Theta(\mbf{x})\right)^T \left(\mbf{x}(\alpha) - \mbf{x} \right) \, .
\end{split}
\end{displaymath}
and thus the inequality (\ref{AR_gradient}) can be rewritten as
\begin{displaymath}
\left(\nabla\Theta(\boldsymbol{\xi}_{\alpha}) - \nabla\Theta(\mbf{x})\right)^T \left(\mbf{x}(\alpha) - \mbf{x} \right) \leqslant (\sigma_G -1) \nabla\Theta(\mbf{x})^T \left(\mbf{x}(\alpha) - \mbf{x} \right)
\end{displaymath}
Since $\sigma_G<1$, from Proposition \ref{prop:proprieties_P} (c) it follows that
\begin{displaymath}
(\sigma_G - 1) \nabla\Theta(\mbf{x})^T \left(\mbf{x}(\alpha) - \mbf{x}\right) =   (1-\sigma_G) \frac{\norm{\mbf{x}(\alpha) - \mbf{x}}^2}{\alpha} 
\geqslant  (1-\sigma_G) \frac{\norm{\mbf{x}(\widetilde{\alpha}) - \mbf{x}}}{\widetilde{\alpha}}\ \norm{\mbf{x}(\alpha) - \mbf{x}} > 0 \, .
\end{displaymath}
The theorem is proved if we show that it exists $\overline{\alpha} \in (0, \widetilde{\alpha}]$ such that for all $\alpha \in [0, \overline{\alpha}]$
\begin{displaymath}
\begin{split}
    \left(\nabla\Theta(\boldsymbol{\xi}_{\alpha}) - \nabla\Theta(\mbf{x})\right)^T \left(\mbf{x}(\alpha) - \mbf{x} \right) \leqslant (1-\sigma_G) \frac{\norm{\mbf{x}(\widetilde{\alpha}) - \mbf{x}}}{\widetilde{\alpha}}\ \norm{\mbf{x}(\alpha) - \mbf{x}} \, .
\end{split}
\end{displaymath}
This follows from the fact that
\begin{displaymath}
\lim_{\alpha \rightarrow 0} \left| \left(\nabla\Theta(\boldsymbol{\xi}_{\alpha}) - \nabla\Theta(\mbf{x})\right)^T \frac{\left(\mbf{x} - \mbf{x}(\alpha) \right)}{\norm{\mbf{x} - \mbf{x}(\alpha)}}\right| \leqslant \lim_{\alpha \rightarrow 0} \norm{\nabla\Theta(\boldsymbol{\xi}_{\alpha}) - \nabla\Theta(\mbf{x})} = 0
\end{displaymath}
where the last equality is a consequence of Proposition \ref{prop:proprieties_P} (b) and of the regularity assumptions on $\Theta$.
\end{proof}
\qed

\begin{theorem}\label{thm:exist_alpha}
Given $\Theta : \mathbb{R}^n \rightarrow \mathbb{R}$, a continuously differentiable function on $\Omega$, $\mbf{x} \in \Omega$, and $\rho \in (0, 1]$, it exists $\overline{\alpha}>0$ so that, for all $\alpha \in (0, \overline{\alpha}]$, $\mathcal{N}_{\alpha} (\mbf{x})= \emptyset$ and thus
\begin{equation}\label{eq:new_cond}
    \sqrt{ \sum_{i \in \mathcal{M}_{\alpha}(\mbf{x})} (P_i(x_i - \partial_i \Theta(\mbf{x})) - x_i)^2  }\geqslant\ \rho \sqrt{ \sum_{i \in \mathcal{N}_{\alpha}(\mbf{x})} (P_i(x_i - \partial_i \Theta(\mbf{x})) - x_i)^2 }  \, .
\end{equation}
\end{theorem}

\begin{proof}
For all $i \in \left\{ 1, \dots, n \right\} \smallsetminus \mathcal{B}(\mbf{x})$ we only have three possibilities:
\begin{itemize} 
    \item $x_i \in \mathring{\Omega_i}$, where $\mathring{\Omega_i}$ denotes the interior of $\Omega_i$. Then, since $\mathring{\Omega}_i$ is an open set, $\exists \, \overline{\alpha}_i > 0$ such that $x_i- \alpha \partial_i \Theta (\mbf{x}) \in \mathring{\Omega}_i \subseteq \Omega_i$ $\forall \alpha \leq \overline{\alpha}_i$. 
    \item $x_i = \ell_i$ and $\partial_i \Theta (\mbf{x}) <0 $. Then $x_i- \alpha \partial_i \Theta (\mbf{x}) \in  \Omega_i$ $\forall \alpha \leq \overline{\alpha}_i := -\frac{u_i - \ell_i}{\partial_i \Theta (\mbf{x})}$.
    \item $x_i = u_i$ and $\partial_i \Theta (\mbf{x}) >0 $. Then $x_i- \alpha \partial_i \Theta (\mbf{x}) \in  \Omega_i$ $\forall \alpha \leq \overline{\alpha}_i := \frac{u_i - \ell_i}{\partial_i \Theta (\mbf{x})}$
    \end{itemize}
Therefore, for all $i \in \left\{ 1, \dots, n \right\} \smallsetminus \mathcal{B}(\mbf{x})$ it exists $\overline{\alpha}_i> 0$ such that $i \in \mathcal{M}_{\alpha}(\mbf{x})$ $\forall \alpha \in (0, \overline{\alpha}_i]$. By choosing $ \overline{\alpha} =\min\limits_{i \in \{1, \dots, n\} \smallsetminus \mathcal{B}(\mbf{x})} \overline{\alpha}_i$, it follows that, for all $\alpha \leq \overline{\alpha}$, $\mathcal{N}_{\alpha} (\mbf{x})= \emptyset$ 
 and thus 
\begin{equation}
\sqrt{ \sum_{i \in \mathcal{M}_{\alpha}(\mbf{x})} (P_i(x_i - \partial_i \Theta(\mbf{x})) - x_i)^2 }  \geqslant 0 =\ \rho \sqrt{\sum_{i \in \mathcal{N}_{\alpha}(\mbf{x})} (P_i(x_i - \partial_i \Theta(\mbf{x})) - x_i)^2 } \, . 
\end{equation}
    Hence the theorem is proved. \qed
\end{proof}

\begin{figure}[H]
\begin{center}
\includegraphics[width=12cm]{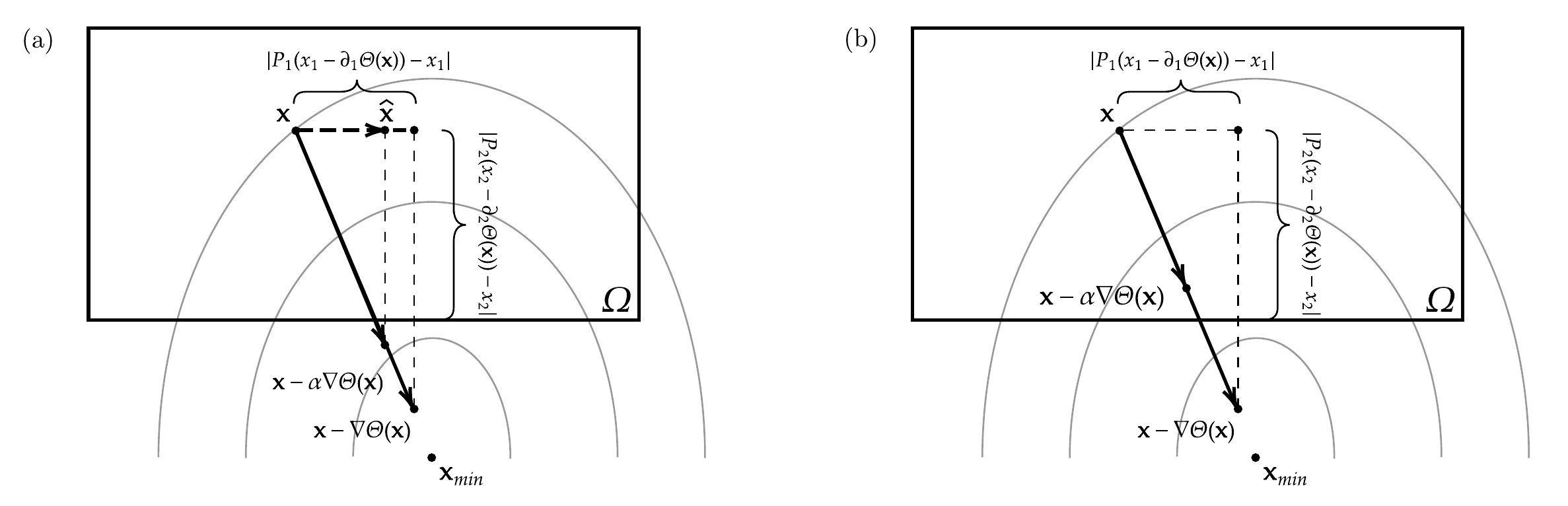}
\end{center}
\caption{Illustration of the benefit of the additional condition (\ref{eq:new_cond}). In (a) only the first component of $\mbf{x}$ is updated as $1 \in \mathcal{M}_{\alpha}(\mbf{x})$ and $2 \in \mathcal{N}_{\alpha}(\mbf{x})$. In this scenario NLPC may get stucked in a point which is not stationary because the chosen stepsize is too big and the second component never updated. As shown in (b), inequality (\ref{eq:new_cond}) prevents this issue by promoting the choice of a smaller stepsize so that an higher number of components is updated. Here, $\widehat{\mbf{x}} = \projP{\mbf{x}-\alpha \nabla \Theta(\mbf{x})}{\mbf{x}}$, $\mbf{x}_{min}$ is a stationary point of $\Theta$, and $\rho=1$. }\label{fig:new_cond} 
\end{figure}

\begin{remark}
The previous theorems hold in particular if $\Theta$ is defined as in equation (\ref{eq:def_theta}). Specifically, inequality (\ref{AR_gradient}) is the classical Armijo rule along the projection arc where we employed the operator introduced in Definition \ref{defP}. Inequality (\ref{eq:new_cond}) is an additional condition that prevents NLPC from choosing a too large stepsize, that would result in an actual update of only few components. An illustrative example can be seen in Figure \ref{fig:new_cond}. 

Theorem \ref{th:GC} and Theorem \ref{thm:exist_alpha} together guarantee that the stepsize within the gradient descent step of the NLPC algorithm is well defined.
\end{remark}

Henceforth, $\left\{ \mbf{x}_k \right\}_{k \in \mathbb{N}} \subseteq \Omega$ and $\left\{ \mbf{\alpha}^{j_k} \right\}_{k \in \mathbb{N}}$ shall denote a sequence of points generated with the NLPC algorithm described in Algorithm \ref{algo:gp_newton}, and the corresponding stepsizes, respectively. In particular, $\alpha \in (0, 1)$, while $j_k$ is a suitable exponent whose value belongs to a different range depending on whether the Newton's or the gradient descent approach has been used at the $k$-th iteration. 

\begin{lemma}\label{red:dis_conds}
Let $\left\{ \mbf{x}_k \right\}_{k \in \mathbb{N}}$ be a sequence generated with the NLPC algorithm. For each $k \in \mathbb{N}$
\begin{description}
\item[(a)] if $\mbf{x}_{k+1}$ has been obtained with a projected gradient descent step, then 
$\Theta(\mbf{x}_{k+1}) \leq \Theta(\mbf{x}_{k})$
\item[(b)] if $\mbf{x}_{k+1}$ has been obtained with a projected Newton's step, then $$\Theta(\mbf{x}_{k+1}) \leq \left( 1-\alpha^J \sigma_N \right)^{n_{k}+1} \Theta(\mbf{x}_0) \, ,$$ 
being $n_{k}$ the number of projected Newton's steps performed until iteration $k$. 
\end{description}
\end{lemma}
\begin{proof}
(a)  By the Armijo rule along the projected arc  in equation (\ref{AR_gradient})
\begin{displaymath}
\Theta(\mbf{x}_{k+1}) \leq \Theta(\mbf{x}_{k}) + \sigma_G \nabla \Theta(\mbf{x}_{k})^T (\mbf{x}_{k+1} - \mbf{x}_{k}) \leq \Theta(\mbf{x}_{k})
\end{displaymath}
where the last inequality follows from Theorem \ref{th:descent_dir}. Thus we have the thesis.\\ 
(b) If $\mbf{x}_{k+1}$ is defined with the projected Newton's method, then it exists $j \in \{ 0, \dots, J \}$ such that
\begin{equation}\label{eq:dis_newton}
\Theta(\mbf{x}_{k+1}) \leq (1-\alpha^j \sigma_N)\ \Theta(\mbf{x}_{k}) \leq (1-\alpha^J \sigma_N)\ \Theta(\mbf{x}_{k}) \, ,
\end{equation}
where in last inequality we exploited the fact that $\alpha < 1$. The thesis follows by iteratively applying (\ref{eq:dis_newton}) for each projected Newton's step, and the results of point (a) for each projected gradient descent step.\\
\qed
\end{proof}

\begin{theorem}\label{th:conv_analysis}
Let $\left\{ \mbf{x}_k \right\}_{k \in \mathbb{N}}$ be a sequence generated with the NLPC algorithm, and let $\mbf{x}^*$ be an accumulation point of $\left\{ \mbf{x}_k \right\}_{k \in \mathbb{N}}$; then $\mbf{x}^*$ is a stationary point of $\Theta$ in $\Omega$. Additionally if the projected Newton's method has been used for infinitely many $k$, then $\mbf{x}^*$ is a solution of (\ref{eq:box_eqs}).
\end{theorem}
\begin{proof}
Since $\mbf{x}^*$ is an accumulation point of $\left\{ \mbf{x}_k \right\}_{k \in \mathbb{N}}$, there exists a subsequence $\left\{ \mbf{x}_{k} \right\}_{k\in K} \subseteq \left\{ \mbf{x}_k \right\}_{k \in \mathbb{N}} $, $K \subseteq \mathbb{N}$, such that $\lim\limits_{k(\in K) \to \infty} \mbf{x}_{k} = \mbf{x}^*$ and thus $\lim\limits_{k(\in K) \to \infty} \Theta(\mbf{x}_{k}) = \Theta(\mbf{x}^*)$.\\
For all $k \in K$, from Lemma $\ref{red:dis_conds}$ it follows
\[
\Theta(\mbf{x}_{k}) \leq \left( 
1-\alpha^J \sigma_N \right)^{n_{k}} \Theta(\mbf{x}_0) \, .
\]
If the projected Newton's method has been used for infinitely many $k$, then $\lim\limits_{k(\in K) \to \infty} n_{k} = + \infty$; therefore
\[
\Theta(\mbf{x}^*) = \lim\limits_{k(\in K) \to \infty} \Theta(\mbf{x}_{k}) \leq  \lim\limits_{k(\in K) \to \infty} \left( 1-\alpha^J \sigma_N \right)^{n_{k}} \Theta(\mbf{x}_0) = 0 \, .
\]
Hence $\Theta(\mbf{x}^*) = 0$, that is $\mbf{x}^*$ solves (\ref{eq:box_eqs}) and is a stationary point of $\Theta$. 

Instead, if the projected gradient direction has been used for all but finitely many iterations, then it exists $ \overline{k} \in \mathbb{N}$ so that $\mbf{x}_{k+1}$ has been obtained through a gradient descent step $\forall\ k \geqslant \overline{k}$. From Lemma \ref{red:dis_conds} it follows $0 \leqslant \Theta(\mbf{x}_{k+1}) \leqslant \Theta(\mbf{x}_{k}) \ \forall \ k \geqslant \overline{k}$, that is  $\{\Theta(\mbf{x}_{k})\}_{k \geqslant \overline{k}}$ is not increasing and bounded below by zero. Hence it converges and 
$$\lim\limits_{k \to \infty} \left( \Theta(\mbf{x}_{k+1}) - \Theta(\mbf{x}_{k}) \right) = 0 \, .$$
Henceforth we shall denote with $\{ \alpha^{j_k} \}_{k \in \mathbb{N}}$ the sequence of stepsizes used within NLPC. From (\ref{AR_gradient}), (\ref{eq9}), and (\ref{eq:prop_norm}) it holds
\begin{displaymath}
\begin{split}
 \Theta(\mbf{x}_{k+1}) - \Theta(\mbf{x}_{k}) & \leqslant \sigma_G \nabla \Theta(\mbf{x}_k)^T \left(\projP{\mbf{x}_k - \alpha^{j_k} \nabla \Theta(\mbf{x}_k)}{\mbf{x}_k} - \mbf{x}_k \right) \\
 & = - \sigma_G \frac{\norm{  \projP{\mbf{x}_k-\alpha^{j_k}\nabla\Theta(\mbf{x}_k)}{\mbf{x}_k}-\mbf{x}_k}^2}{\alpha^{j_k}} \\
 & = - \sigma_G\ \alpha^{j_k} \sum_{i \in \mathcal{M}_{\alpha^{j_k}}(\mbf{x}_k)} (\partial_i \Theta(\mbf{x}_k))^2 \leq 0 
\end{split}
\end{displaymath}
and thus 
\begin{equation}\label{eq:dim_limit}
\lim\limits_{k (\in K) \to \infty}
 \alpha^{j_k} \sum_{i \in \mathcal{M}_{\alpha^{j_k}}(\mbf{x}_k)} (\partial_i \Theta(\mbf{x}_k))^2 = 0 \, .
\end{equation} 
Two cases exist:  $\liminf\limits_{k(\in K) \to \infty} \alpha^{j_k} > 0$ (\textit{case 1}) and $\liminf\limits_{k(\in K) \to \infty} \alpha^{j_k} = 0$ (\textit{case 2}).

If \textit{case 1} holds, then equation (\ref{eq:dim_limit}) implies
\begin{displaymath}
\begin{split}
0 & = 
\lim\limits_{k(\in K) \to \infty} \sum_{i \in \mathcal{M}_{\alpha^{j_k}}(\mbf{x}_{k})} (\partial_i \Theta(\mbf{x}_{k}))^2 \\
& \geq \lim\limits_{k(\in K) \to \infty} \sum_{i \in \mathcal{M}_{\alpha^{j_k}}(\mbf{x}_{k})} (P_i(x_{{k},i} - \partial_i \Theta(\mbf{x}_{k})) - x_{{k},i})^2 \\
& \geq \rho \lim\limits_{k(\in K) \to \infty} \sum_{i \in \mathcal{N}_{\alpha^{j_k}}(\mbf{x}_{k})} (P_i(x_{{k},i} - \partial_i \Theta(\mbf{x}_{k})) - x_{{k},i})^2 \, ,
\end{split}
\end{displaymath}
where the last inequality comes from the constraint described by (\ref{eq:new_cond}). Hence, in particular
\begin{displaymath}
\begin{split}
\lim\limits_{k(\in K) \to \infty} & \sum_{i \in \mathcal{M}_{\alpha^{j_k}}(\mbf{x}_{k})} (P_i(x_{{k},i} - \partial_i \Theta(\mbf{x}_{k})) - x_{{k},i})^2  = \\
& \lim\limits_{k(\in K) \to \infty} \sum_{i \in \mathcal{N}_{\alpha^{j_k}}(\mbf{x}_{k})} (P_i(x_{{k},i} - \partial_i \Theta(\mbf{x}_{k})) - x_{{k},i})^2 = 0 \, .
\end{split}
\end{displaymath}

Since additionally $(P_i(x_{{k},i} - \partial_i \Theta(\mbf{x}_{k})) - x_{{k},i}) = 0$ for all $i \in \mathcal{B}(\mbf{x}_{k})$, from equation (\ref{sets_BMN}), from the continuity of $\nabla \Theta$ and of the classical projector, and being $\mbf{x}^*$ the limit point of $\{\mbf{x}_{k}\}_{k \in K}$, it follows
\begin{displaymath}
||P(\mbf{x}^* - \nabla \Theta(\mbf{x}^*)) - \mbf{x}^*||^2  = \lim\limits_{k(\in K) \to \infty} \sum_{i=1}^n (P_i(x_{{k},i} - \partial_i \Theta(\mbf{x}_{k})) - x_{{k},i})^2 = 0 \, .
\end{displaymath}
Hence $\mbf{x}^*$ is a stationary point.

On the other hand, \textit{case 2} implies it exists an infinite set $K' \subset K$ such that
$\lim\limits_{k(\in K') \to \infty} \alpha^{j_k} = 0$, and thus $\lim\limits_{k(\in K') \to \infty} \alpha^{j_k-1} = 0$. Therefore, by defining $J = \left\{ i \in \{1, \dots, n\}\ s.t.\ i \notin \mathcal{B}(\mbf{x}^*)\ \land \  |\partial_i \Theta(\mbf{x}^*)| > 0 \right\}$, the following holds.
\begin{description}
    \item[(i)] $\mbf{x}^*$ is a stationary point iff $J = \emptyset$. \\
    More specifically, $|P_i(x^*_{i} - \partial_i \Theta(\mbf{x}^*)) - x^*_i| = 0$ iff $i \notin J$.
    \item[(ii)] It exists $\overline{k}$ such that $\forall\ k \in K',\ k \geq \overline{k}$, $J \subseteq \mathcal{M}_{\alpha^{j_k-1}}(\mbf{x}_k)$. \\ Indeed, let's consider $i \in J$. Since in particular $i \notin \mathcal{B}(\mbf{x}^*)$ and $\alpha<1$, from Theorem \ref{thm:exist_alpha} it follows that it exists $\overline{j} \in \mathbb{N}$ such that $\mathcal{N}_{\alpha^j}(\mbf{x}^*) = \emptyset$ and $i \in \mathcal{M}_{\alpha^j}(\mbf{x}^*)$ $\forall j \geqslant \overline{j}$.
    It can be easily shown that, being $\mbf{x}^*$ the limit point of $\left\{ \mbf{x}_{k} \right\}_{k\in K}$, this implies it exists $k_i$ such that $\forall\ k \in K'$ with $\ k \geq k_i$, $i \in \mathcal{M}_{\alpha^{j}}(\mbf{x}_k)$ $\forall j > \overline{j}$. Additionally, since $\lim\limits_{k(\in K') \to \infty} \alpha^{j_k-1} = 0$, it exists $k'_i \geq k_i $ such that $\forall\ k \in K',\ k \geq k'_i$, $\alpha^{j_k-1} < \alpha^{\overline{j}}$. Hence the thesis follows by considering $\overline{k} = \max\limits_{i \in \{1, \dots, n\}}\{ k'_i\}$.
\end{description}
To prove that $\mbf{x}^*$ is a stationary point we proceed by contradiction and we assume that it exists $i \in J$. Then, by the results in (i) and (ii), it follows
\begin{displaymath}
\begin{split}
\lim\limits_{k(\in K') \to \infty} & \sqrt {\sum_{i \in \mathcal{M}_{\alpha^{j_k-1}}(\mbf{x}_{k})} |P_i(x_{k,i}-\partial_i \Theta(\mbf{x}_{k}))-x_{k,i}|^2 } \\ 
& \geq \lim\limits_{k(\in K') \to \infty} \sqrt{ \sum_{i \in J} |P_i(x_{k,i}-\partial_i \Theta(\mbf{x}_{k}))-x_{k,i}|^2 } \\
& =  \sqrt{ \sum_{i \in J} |P_i(x^*_{i}-\partial_i \Theta(\mbf{x}^*))-x^*_{i}|^2 } > 0 \, ,
\end{split}
\end{displaymath}
while, denoted $J^C = \left\{1, \dots, n \right\} \setminus J$,
\begin{displaymath}
\begin{split}
\lim\limits_{k(\in K') \to \infty} & \rho \sqrt{ \sum_{i \in \mathcal{N}_{\alpha^{j_k-1}}(\mbf{x}_{k})} |P_i(x_{k,i}-\partial_i \Theta(\mbf{x}_{k}))-x_{k,i}|^2 } \\ 
& \leq \rho \lim\limits_{k(\in K') \to \infty} \sqrt{ \sum_{i \in J^C} |P_i(x_{k,i}-\partial_i \Theta(\mbf{x}_{k}))-x_{k,i}|^2 } \\
& = \rho \sqrt{ \sum_{i \in J^C} |P_i(x^*_{i}-\partial_i \Theta(\mbf{x}^*))-x^*_{i}|^2 } = 0 \, .
\end{split}
\end{displaymath}
Therefore, for sufficiently large $k \in K'$
\begin{displaymath}
\sqrt{ \sum_{i \in \mathcal{M}_{\alpha^{j_k-1}}(\mbf{x}_{k})} |P_i(x_{k,i}-\partial_i \Theta(\mbf{x}_{k}))-x_{k,i}|^2 } \geq\ \rho \sqrt{ \sum_{i \in \mathcal{N}_{\alpha^{j_k-1}}(\mbf{x}_{k})} |P_i(x_{k,i}-\partial_i \Theta(\mbf{x}_{k}))-x_{k,i}|^2 }
\end{displaymath}
i.e. condition (\ref{eq:new_cond}) is satisfied by the stepsize $\alpha^{j_k-1}$ which is the last stepsize tried by NLPC before the chosen one. As a consequence, such a stepsize cannot satisfy condition (\ref{AR_gradient}), i.e.
\begin{displaymath}
\begin{split}
\Theta( \projP{\mbf{x}_k & -\alpha^{j_k-1} \nabla \Theta(\mbf{x}_k)}{\mbf{x}_k})  - \Theta(\mbf{x}_k)\\
& > \sigma_G \nabla \Theta(\mbf{x}_k)^T \left(\projP{\mbf{x}_k - \alpha^{j_k-1} \nabla \Theta(\mbf{x}_k)}{\mbf{x}_k} - \mbf{x}_k \right) \, .
\end{split}
\end{displaymath}
By the mean value theorem, it exists $\tau \in (0, 1)$ such that, defined $\boldsymbol{\xi}_k =\tau \mbf{x}_{k} + (1-\tau)  \projP{\mbf{x}_{k} - \alpha^{j_k-1} \nabla \Theta(\mbf{x}_{k})}{\mbf{x}_{k}}$, then 
\begin{displaymath}
\begin{split}
\Theta( \projP{\mbf{x}_k & -\alpha^{j_k-1} \nabla \Theta(\mbf{x}_k)}{\mbf{x}_k})   - \Theta(\mbf{x}_k) \\
& = \nabla \Theta(\boldsymbol{\xi}_k)^T \left(\projP{\mbf{x}_k - \alpha^{j_k-1} \nabla \Theta(\mbf{x}_k)}{\mbf{x}_k} - \mbf{x}_k \right) \\
& = \left(\nabla \Theta(\boldsymbol{\xi}_k) - \nabla \Theta(\mbf{x}_k)\right)^T \left(\projP{\mbf{x}_k - \alpha^{j_k-1} \nabla \Theta(\mbf{x}_k)}{\mbf{x}_k} - \mbf{x}_k \right) \\
& + \nabla \Theta(\mbf{x}_k)^T \left(\projP{\mbf{x}_k - \alpha^{j_k-1} \nabla \Theta(\mbf{x}_k)}{\mbf{x}_k} - \mbf{x}_k \right) \, .
\end{split}
\end{displaymath}
Together with the previous result this implies
\begin{displaymath}
\begin{split}
(1-\sigma_G) & \nabla \Theta(\mbf{x}_k)^T \left(\mbf{x}_k - \projP{\mbf{x}_k - \alpha^{j_k-1} \nabla \Theta(\mbf{x}_k)}{\mbf{x}_k} \right) \\
& < \left(\nabla \Theta(\boldsymbol{\xi}_k) - \nabla \Theta(\mbf{x}_k)\right)^T \left(\projP{\mbf{x}_k - \alpha^{j_k-1} \nabla \Theta(\mbf{x}_k)}{\mbf{x}_k} - \mbf{x}_k \right) \\
& \leq \norm{\nabla \Theta(\boldsymbol{\xi}_k) - \nabla \Theta(\mbf{x}_k)} \cdot \norm{ \projP{\mbf{x}_k - \alpha^{j_k-1} \nabla \Theta(\mbf{x}_k)}{\mbf{x}_k} - \mbf{x}_k} \, , 
\end{split}
\end{displaymath}
hence 
\begin{displaymath}
\begin{split}
\frac{1}{1-\sigma_G} \norm{\nabla \Theta(\boldsymbol{\xi}_k) - \nabla \Theta(\mbf{x}_k)} 
& > \frac{\nabla \Theta(\mbf{x}_k)^T \left(\mbf{x}_k - \projP{\mbf{x}_k - \alpha^{j_k-1} \nabla \Theta(\mbf{x}_k)}{\mbf{x}_k} \right) }{\norm{ \projP{\mbf{x}_k - \alpha^{j_k-1} \nabla \Theta(\mbf{x}_k)}{\mbf{x}_k} - \mbf{x}_k}}\\
& =  \sqrt{\sum_{i \in \mathcal{M}_{\alpha^{j_k-1}}(\mbf{x}_{k})} (\partial_i \Theta (\mbf{x}_{k}))^2 } \, , 
\end{split}
\end{displaymath}
where the last equality comes from equations (\ref{eq9}) and (\ref{eq:prop_norm}).
Therefore, from the properties of the classical projector and from the result previously shown in (ii), it follows
\begin{displaymath}
\begin{split}
    0 & = \lim\limits_{k(\in K') \to \infty}  \sqrt{\sum_{i \in \mathcal{M}_{\alpha^{j_k-1}}(\mbf{x}_{k})} (\partial_i \Theta (\mbf{x}_{k}))^2 } \\
    & \geq \lim\limits_{k(\in K') \to \infty}  \sqrt{\sum_{i \in \mathcal{M}_{\alpha^{j_k-1}}(\mbf{x}_{k})} |P_i(x_{k,i}-\partial_i \Theta(\mbf{x}_{k}))-x_{k,i}|^2 } \\
    & \geq \sqrt{\sum_{i \in J} |P_i(x^*_{i}-\partial_i \Theta(\mbf{x}^*))-x^*_{i}|^2} \, ,
\end{split}
\end{displaymath}

which is possible only if $J=\emptyset$ and thus contradicts our hypothesis.
\qed
\end{proof}

As a final remark, we observe that all the results proven in this section can be easily extended to the case where the gradient direction is normalized \cite{Wattetal2020}.

\section{Application to chemical reaction networks}\label{sec:NLPC_for_CRN}
Let's consider a chemical reaction network (CRN) composed of $r$ chemical reactions involving $n$ well-mixed proteins. Specifically, in this work we will focus on the CRN devised for modeling cell signaling during the G1-S transition phase in colorectal cells described in \cite{Tortolina2015,sommariva2021_scirep} and henceforth denoted as CR-CRN. In this case $n=419$ and $r=851$.

By assuming that the law of mass action holds \cite{Yu_Craciun_2018,sommariva2021_JMB}, the dynamics of the CRN gives rise to a set of $n$ ordinary differential equations (ODEs)
\begin{equation}
\dot{\mbf{x}} = \mbf{S} \mbf{v}(\mbf{x}, \mbf{k})\label{eq:ODEs_CRN}
\end{equation}
where the state vector $\mbf{x} \in \mathbb{R}^n_+$ contains the protein molecular concentrations (nM); the superposed dot denotes the time derivative; $\mbf{S}$ is the constant stoichiometric matrix of size $n \times r$; $\mbf{k} \in \mathbb{R}_+^r$ are the rate constants of the reactions; and $\mbf{v}(\mbf{x}, \mbf{k}) \in \mathbb{R}^r_+$ is the time-variant vector of the reaction fluxes. Specifically, from the law of mass action it follows \cite{otero2017}
\begin{equation}\label{eq:def_v}
\mbf{v}(\mbf{x}, \mbf{k}) = \textrm{diag}(\mbf{k}) \mbf{z}(\mbf{x})
\end{equation}
where the elements of $\mbf{z}(\mbf{x})$ are monomials of the form $z_j(\mbf{x}) = \prod_{i=1}^n x_i^{p_{ij}}$, $\forall j = 1, \dots, r$. In the CR-CRN, $p_{ij} \in \left\{0, 1, 2 \right\}$, because all the reactions involve up to two reactants. 

Given a solution $\mbf{x}(t)$ of system (\ref{eq:ODEs_CRN}), a semi-positive conservation vector is a constant vector $\boldsymbol{\gamma} \in \mathbb{N}^n \setminus \{\mathbf{0}\}$ for which it exists $c \in \mathbb{R}_+$ so that $\boldsymbol{\gamma}^T \mathbf{x}(t) = c$ $\forall \ t$ \cite{sommariva2021_JMB,shinar2009}. Conservation vectors can be determined by studying the kernel of $\mathbf{S}^T$ \cite{schuster1991}. In the remaining of the paper we shall assume that the considered CRN satisfies the following properties in terms of its conservation vectors.
\begin{description}
\item{(i)} The CRN is \textit{weakly elemented} \cite{sommariva2021_JMB}, i.e. it exists a set of independent generators $\left\{ \boldsymbol{\gamma}_1, \dots, \boldsymbol{\gamma}_p \right\} \subset \mathbb{N}^n \setminus \{\mbf{0}\}$ of the semi-positive conservation vectors such that $p = n - \text{rank}(\mathbf{S})$ and, up to a change  of the proteins order, 
\begin{equation}
    \mathbf{N} :=
 \begin{bmatrix} 
  \boldsymbol{\gamma}_1^T \\
  \vdots \\
  \boldsymbol{\gamma}_{p}^T \end{bmatrix} 
  = \left[\mathbf{I}_p, \mathbf{N}_2 \right] \ , 
\end{equation}
being $\mathbf{I}_p$ the identity matrix of size $p \times p$.
\item{(ii)} The CRN satisfies the \textit{global stability condition} \cite{sommariva2021_JMB}, i.e. for each $\mbf{c} \in \mathbb{R}^p_+$ it exists a unique asymptotically stable state on the stoichiometric compatibility class (SCC) $\left\{\mbf{x} \in \mathbb{R}^n_+\ \text{s.t.}\ \mbf{N}\mbf{x} = \mbf{c} \right\}$. Fixed a SCC, the corresponding asymptotically stable state $\mbf{x}_e \in \mathbb{R}^n_+$ solves the system
\begin{equation}\label{eq:rec_system}
\begin{cases}
\mbf{S} \mbf{v}(\mbf{x}, \mbf{k}) = 0 \\
\mbf{N} \mbf{x} - \mbf{c} = 0.
\end{cases}
\end{equation}
\end{description}

\begin{lemma}\label{lemma:squared_system}
For a weakly elemented CRN satisfying the global stability condition, the system in (\ref{eq:rec_system}) is equivalent to the square system 
\begin{equation}\label{eq:square_system}
\begin{cases}
\mbf{S}_2 \mbf{v}(\mbf{x}, \mbf{k}) = 0 \\
\mbf{N} \mbf{x} - \mbf{c} = 0
\end{cases}
\end{equation}
where $\mbf{S}_2$ is a matrix of size $(n-p) \times r$ defined by the last $n-p$ rows of $\mbf{S}$.
\end{lemma}
\begin{proof}
Obviously, a solution of (\ref{eq:rec_system}) also solves (\ref{eq:square_system}). On the other hand, let $\mbf{x}_e$ be a solution of (\ref{eq:square_system}). The theorem is proved by showing that 
\begin{displaymath}
\mbf{S}_1  \mbf{v}(\mbf{x}_e, \mbf{k}) = 0 \ ,  
\end{displaymath}
being $\mbf{S}_1$ the matrix of size $p \times r$ defined by the first $p$ rows of $\mbf{S}$. To this end we observe that, since any conservation vector belongs to the kernel of $\mathbf{S}^T$, it holds 
\begin{displaymath}
\mbf{0} = \mbf{N} \mbf{S} = \mbf{S}_1 + \mbf{N_2} \mbf{S_2} \ .
\end{displaymath}
Therefore 
\begin{displaymath}
\mbf{S}_1  \mbf{v}(\mbf{x}_e, \mbf{k}) = - \mbf{N_2} \mbf{S_2} \mbf{v}(\mbf{x}_e, \mbf{k}) = 0\ .  
\end{displaymath}
\qed
\end{proof}
According to Lemma \ref{lemma:squared_system}, in a weakly elemented CRN satisfying the global stability condition, the equilibrium point on a fixed SCC can be computed by solving a box--constrained system as in equation (\ref{eq:box_eqs}), being $\Omega = \mathbb{R}_+^n$ and 
\begin{equation}\label{eq:f_CRN}
    \mbf{f}(\mbf{x}) = \left[
    \begin{array}{c}
         \mbf{S}_2 \mbf{v}(\mbf{x}, \mbf{k})  \\
         \mbf{N} \mbf{x} - \mbf{c}
    \end{array}\right] \, .
\end{equation}
\begin{lemma} Consider the function $\mbf{f}:\mathbb{R}^n \rightarrow \mathbb{R}^n$ defined as in equation (\ref{eq:f_CRN}). $\mbf{f}$ is continuously differentiable on $\mathbb{R}^n_+$ and
\begin{equation}\label{eq:J_f_CRN}
    \mbf{J}_{\mbf{f}}(\mbf{x}) = \left[
    \begin{array}{c}
         \mbf{S}_2 \textrm{diag}(\mbf{k}) \mbf{J}_\mbf{z}(\mbf{x})  \\
         \mbf{N}  
    \end{array}\right],
\end{equation}
where $[J_{\mbf{z}}(\mbf{x})]_{ji} = p_{ij} x_i^{p_{ij}-1} \prod_{\ell=1, \ell \neq i}^n x_\ell^{p_{\ell j}}$, $\forall i \in \{ 1, \dots, n\}$ and $j = 1, \dots, r$.
\end{lemma}
\begin{proof}
The thesis follows from the definition of the reaction fluxes in
 Eq. (\ref{eq:def_v}). \qed
\end{proof}

\section{Numerical results on the CR-CRN}\label{sec:results}

\subsection{General consideration}\label{sec:implementation}
To show the advantages of using NLPC for computing the asymptotically stable states of a CRN, we applied it to the CR-CRN. The parameters describing the network in a physiological state have been extensively described in previous works \cite{Tortolina2015,sommariva2021_scirep,sommariva2021_JMB} and can be downloaded from the GitHub repository \url{https://​github.​com/​theMI​DAgro​up/​CRC_​CRN.git} as MATLAB\textsuperscript{\textregistered} structure. This includes the list of proteins and reactions involved in the network, as well as the values of the rate constants $\mathbf{k}$ and of the total conserved moieties $\mathbf{c}$. The corresponding stoichiometric matrix $\mathbf{S}$ and reaction fluxes $\mathbf{v}(\mathbf{x}, \mathbf{k})$ can be derived as described in the previous section. The aforementioned repository also contains the MATLAB\textsuperscript{\textregistered} codes implementing the NLPC algorithm and the analysis shown in this paper.

We exploited the model introduced by Sommariva and colleagues \cite{sommariva2021_scirep,sommariva2021_JMB} to test the proposed approach under different biologically-plausible conditions. Specifically, we modified the values of the parameters $\mathbf{k}$ and $\mathbf{c}$ as described in \cite{sommariva2021_scirep} to simulate the effect of some of the mutations that most commonly arise in colorectal cancer. A total of 9 different mutations was considered (loss of function of APC, AKT, SMAD4, PTEN, p53 and gain of function of k-Ras, Raf, PI3K, Betacatenin) which give rise to as many different mutated networks.

From a practical point of view, if not otherwise specified, the parameters required in input by Algorithm \ref{algo:gp_newton} were set as follows. The threshold within the stopping criterion was $\tau=10^{-12}$, while $\sigma_N = \sigma_G = 10^{-4}$ and $\rho = 10^{-2}$. The initial stepsize was $\alpha=0.79$ and a maximum of $J=20$ stepsizes was tested within each iteration of the Newton's method. NLPC is initialized with a point $\mathbf{x}_0$ randomly drawn from the SCC $\left\{\mbf{x} \in \mathbb{R}^n_+\ \text{s.t.}\ \mbf{N}\mbf{x} = \mbf{c} \right\}$ by exploiting the procedure presented in \cite{sommariva2021_JMB}. Additionally,  we  only retained points such that the condition number of the Jacobian matrix $\mbf{J}_{\mbf{f}}(\mbf{x}_0)$ was lower than $10^{17}$. To  avoid the algorithm getting stuck in a stationary point that is not a zero of $\mathbf{f}$ we also set a maximum number of allowed iterations: if the stopping criterion was not reached after 250 iterations, then a new initial point $\mbf{x}_0$ was drawn from the SCC and NLPC was restarted.
Finally, to speed up the performance of NLPC, within the gradient descent method we normalized the gradient direction \cite{Wattetal2020} and we fixed a maximum number of tested stepsizes also for this approach: if conditions (\ref{AR_gradient}) and (\ref{eq:new_cond}) were not met after 40 possible values of the step length, we chose the last tested value and at the following NLPC iteration we performed again a gradient descent step.

\subsection{Comparison with a classical dynamic approach}
\label{par:dyn}
A classical approach \cite{sommariva2021_JMB,sommariva2021_scirep} for computing the stationary state of system (\ref{eq:ODEs_CRN}) on a given SCC consists in simulating the whole concentration dynamics $\mbf{x}(t)$ by solving the Cauchy problem
\begin{equation}\label{eq:chauchy_pb}
\begin{cases}
\dot{\mbf{x}} = \mbf{S} \mbf{v}(\mbf{x}, \mbf{k}) \\
\mbf{x}(0) = \mbf{x}_0
\end{cases} \, ,
\end{equation}
where $\mbf{x}_0$ is a point on the SCC, and then computing the asymptotic value 
\begin{equation}\label{eq:dyn_sol}
\mbf{x}_{dyn} = \lim_{\mbf{t} \to +\infty} \mbf{x}(t)\ .
\end{equation}
In this section we compare the results obtained through this approach with those from NLPC algorithm.

To this end, we started from the CR-CRN and we built 10 different experiments, by varying the values of the kinetic parameters $\mbf{k}$ and of the total conserved moieties $\mbf{c}$ that define the SCC, so as to mimic a colorectal cell either healthy or affected by one of the 9 mutations listed in Subsect. \ref{sec:implementation}. For each experiment, we sampled 50 initial points $\mbf{x}_0^{(j)}$ on the corresponding SCC. For each initial point, i.e. for $j = 1, \dots, 50$, we computed the solution $\mbf{x}_{nlpc}^{(j)}$ provided by the NLPC algorithm and we compared it with the asymptotically stable state $\mbf{x}_{dyn}^{(j)}$ computed through the dynamic approach just described. Specifically, as in \cite{sommariva2021_JMB}, we used the MATLAB\textsuperscript{\textregistered} tool \texttt{ode15s} \cite{Shampine} to integrate the ODEs system in (\ref{eq:chauchy_pb}) on the interval $[0, 2.5 \cdot 10^7]$ and we defined $\mbf{x}_{dyn}^{(j)}$ as the value of the computed solution at the last time-point of the interval.

\begin{figure}[ht]
\begin{center}
\includegraphics[width=10cm]{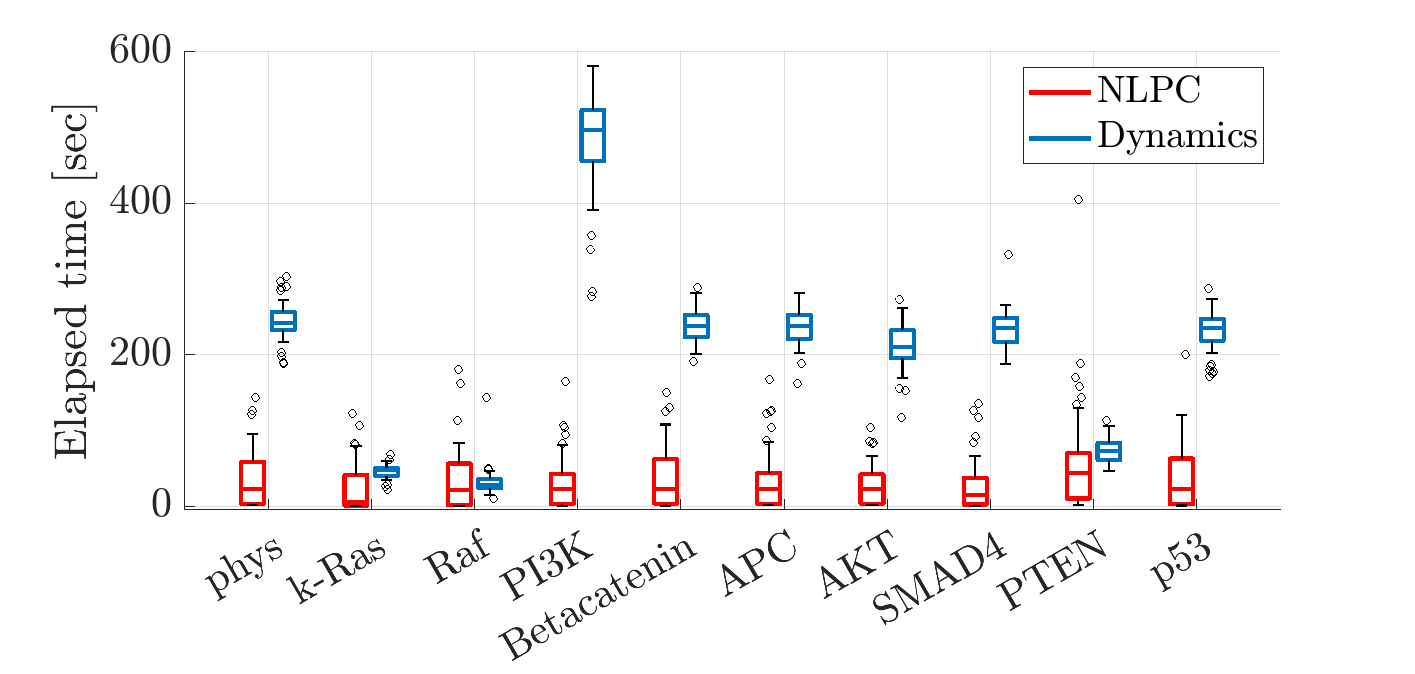}
\end{center}
\caption{Elapsed time for the NLPC algorithm to converge compared to the time required to compute the equilibrium point by solving the dynamical system in (\ref{eq:chauchy_pb}). Boxplots summarize the values obtained across 50 different runs for 10 distinct networks mimicking either a physiological state (phys) or a mutation affecting the protein shown in the axis labels.}\label{fig:time}
\end{figure}

\begin{figure}[ht]
\begin{center}
\includegraphics[width=10cm]{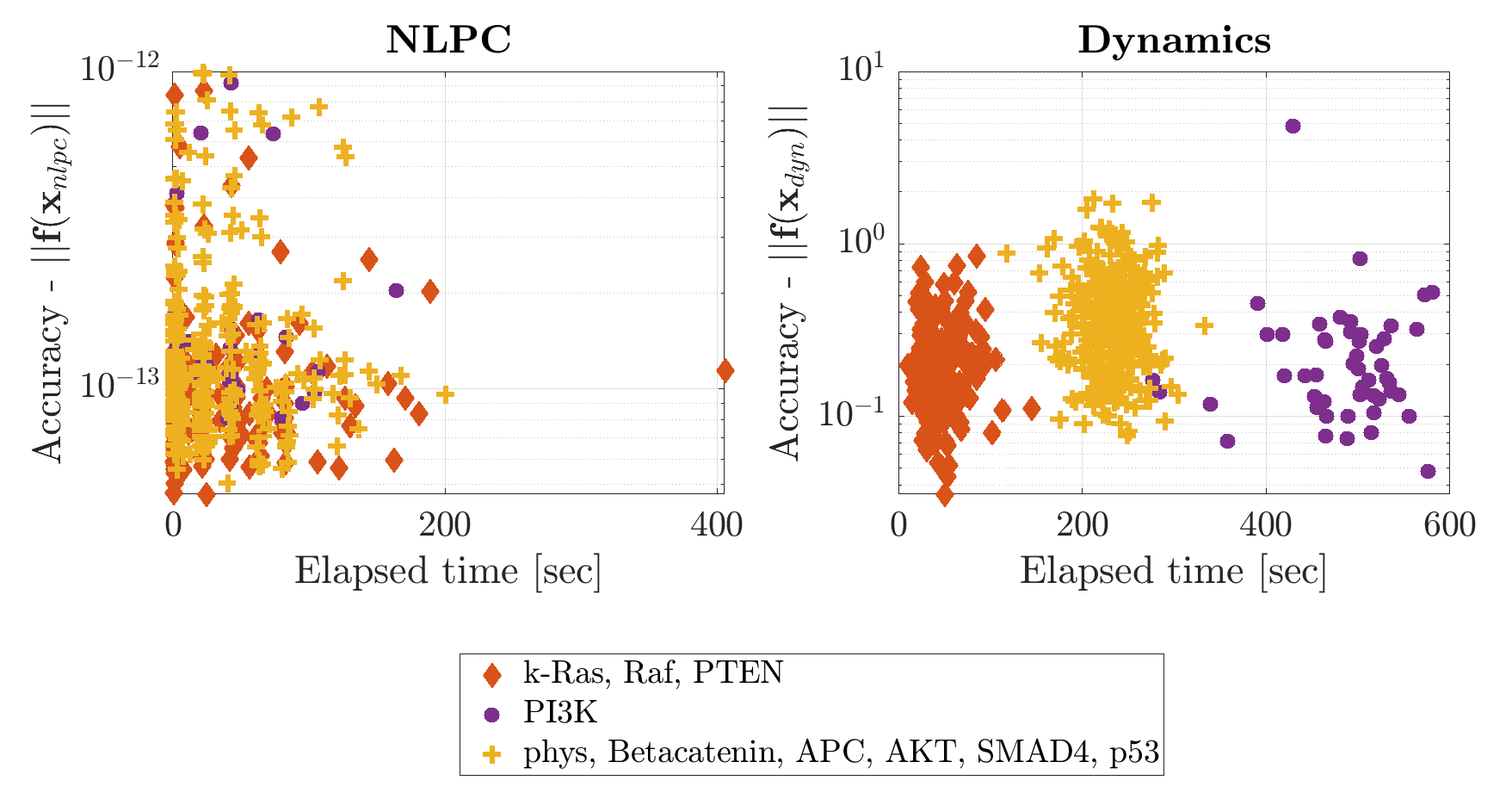}
\end{center}
\caption{Accuracy as a function of the elapsed time for the NLPC algorithm (left) and the dynamic approach (right). Accuracy is quantified as the norm of $\mbf{f}$ evaluated in the results provided by the two algorithms, $\mbf{x}_{nlpc}$ and $\mbf{x}_{dyn}$, respectively. In each panel, 50 different results are shown for each of the considered CRNs that mimic mutation of k-Ras, Raf and PTEN (orange diamonds), physiological state and mutation of Betacatenin, APC, AKT, SMAD4, PTEN, p53 (yellow crosses), and mutation of PI3K (purple dots). This color code as been chosen so as to cluster together results for which the times required for computing $\mbf{x}_{dyn}$ were similar, as depicted in Fig. \ref{fig:time}.
Notice the different scale on the y-axis.}\label{fig:time_precision}
\end{figure}

As shown in Fig. \ref{fig:time} and \ref{fig:time_precision}, NLPC outperforms the dynamic approach in terms of both accuracy of the obtained results and computational cost. Indeed, Fig. \ref{fig:time} shows that in all the 10 considered experiments, the elapsed time for the NLPC algorithm, averaged across 50 runs obtained by varying the initial points, ranges from about $5$ sec (mutated network with gain of function of k-Ras) to $33$ sec (mutated network with loss of function of PTEN). On the contrary, the results of the dynamic approach show an higher variability across the different CRNs and the averaged elapsed time scales up to about $8$ min in the network incorporating a gain of function mutation of PI3K. It is worth noticing that, for each of the 10 experiments, few runs of NLPC required an higher elapsed time (higher than the third quartile of the corresponding distributions). These runs needed a large number of restarts of the NLPC algorithm due to the fact that the maximum number of 250 iterations was reached without meeting the stopping criterion on the norm of $\mbf{f}$, probably because the gradient method tended to stationary points that were not roots of $\mbf{f}$. Future work will be devoted to refining the stopping criterion so that, when needed, NLPC is restarted before reaching 250 iterations.

Since we are looking for the roots of $\mbf{f}$, the accuracy of the obtained results was evaluated by computing the $\ell_2$-norm of $\mbf{f}$ in the solutions provided by the two algorithms, namely $\mbf{x}_{nlpc}^{(j)}$ and $\mbf{x}_{dyn}^{(j)}$ , $j \in \{1, \dots, 50\}$. As shown in Fig. \ref{fig:time_precision}, for all 10 considered experiments the norm of $\mbf{f}$ in the NLPC solutions, $\mbf{x}_{nlpc}^{(j)}$, was always below $10^{-12}$ as imposed by the stopping criterion of the algorithm. Instead the value of $||\mbf{f}(\mbf{x}_{dyn}^{(j)})||$ ranged between $10^{-2}$ to $10^1$, regardless of the time employed to compute the solution $\mbf{x}_{dyn}^{(j)}$.

\subsection{Benefits of the operator $\mathcal{P}$ over the classical projector}
The goal of this section is to quantify the benefit of using the operator $\mathcal{P}$ instead of the classical projector  $P$ on the closed convex set $\Omega$ defined in Eq. (\ref{eq:classical_P}). To this end, for each of the 10 experiments defined in the previous section, and for each of the 50 initial points $\mathbf{x}_0^{(j)}$, $j \in \{1, \dots, 50\}$, drawn on the corresponding SCCs, we computed the solution of NLPC by replacing in Algorithm \ref{algo:gp_newton} the proposed operator $\mathcal{P}$ with the classical projector $P$. We denoted with $\mbf{x}_{ort}^{(j)}$ the corresponding solution.

\begin{figure}[h]
\begin{center}
\includegraphics[width=10cm]{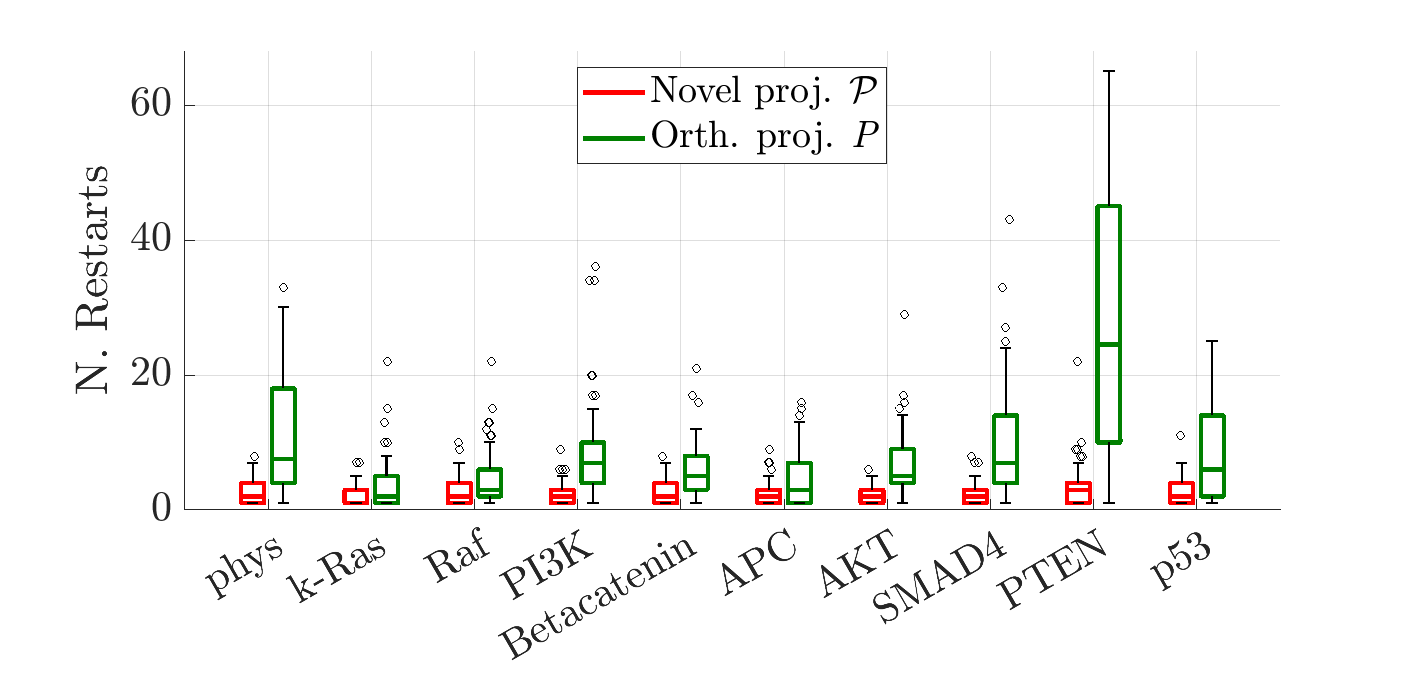}
\end{center}
\caption{Number of restarts required by NLPC in order to satisfy the stopping criterion within the fixed maximum number of iterations. The boxplots describe the values obtained across 50 different runs for 10 distinct networks when, within NLPC, we employed the proposed non linear projector $\mathcal{P}$ (red) or the classical orthogonal projector $P$ (green).}\label{fig:attempts}
\end{figure}

As shown in Figure \ref{fig:attempts}, if combined with the classical projector, NLPC algorithm requires a higher number of restarts and thus a higher elapsed time than those required when the proposed operator is used. Specifically, the ratio between the number of restarts required by the projector $P$ and the one required by the operator $\mathcal{P}$, averaged over all the 10 considered experiments and all the sampled initial points, is around 4.87.  

The bad performances of the projector $P$ are caused by the fact that at any given iteration $k$ all the negative components of the novel proposed point $\mbf{x}_{k+1}$ are set equal to zero. As a consequence, the percentage of proteins estimated as having a null concentration increases sharply and this results in a high condition number of the corresponding Jacobian matrix $\textbf{J}_{\textbf{f}}$ defined as in \eqref{eq:J_f_CRN}.
In turn, the ill-conditioning of $\textbf{J}_{\textbf{f}}$ compromises the stability of Newton's method and thus NLPC algorithm tends to spend most of the allowed iterations by performing gradient descent steps.
As shown in Table \ref{tab:proj_ort}, the use of the operator $\mathcal{P}$ helps preventing this issue.

\begin{table}[ht]
 \centering 
 \caption{Average and standard deviation over 50 initial points of the maximum number of null components (first row) and the maximum condition number of the Jacobian matrix $\mbf{J}_{\mbf{f}}$ (second row) reached across the iterations performed by NLPC. Results obtained by using the novel non liner projector $\mathcal{P}$ (first column) and the classical orthogonal projector $P$ (second column) are compared. Since results across the 10 considered experiments were similar, only those concerning the original physiological CR-CRN are shown. }\label{tab:proj_ort}
	\begin{tabular}{c|*{2}{c|}}	\cline{2-3} 	& \textbf{Novel proj.} $\mathcal{P}$ & \textbf{Orth. proj. $P$} \\ 
\hline	
\multicolumn{1}{|c|}{\textbf{Num. null components (\%)}} 	&0.59 $\pm$ 0.14& 36.99 $\pm$ 8.42 \\\hline
\multicolumn{1}{|c|}{\textbf{Cond. Number $\mbf{J}_{\mbf{f}}$ (log. scale)}} & 14 $\pm$ 2& 17 $\pm$ 3\\  
 \hline  
 \end{tabular}	
  \end{table}

\section{Conclusions}\label{sec:conclusions}
In this paper an iterative algorithm for solving  rootfinding box--constrained problems is presented. It combines both Newton’s and gradient descent methods and exploits the operator $\mathcal{P}$ in Def. \eqref{defP} for assuring the required constraints at each iteration (and preventing numerical instability issues that would occur if the projector $P$ was applied).
Together with a suitable backtracking rule we prove that the method converges to a stationary point of the objective function in Eq. \eqref{eq:def_theta}. Despite outperforming the dynamic approach both in accuracy and speed, in CRNs’ framework the NLPC algorithm provides less information than simulating the whole concentration dynamics. 
However, in many contexts such as tuning kinetic parameters starting from experimental data or for topics described in \cite{sommariva2021_scirep,sommariva2021_JMB} the comprehension of the whole dynamic is not required, but only knowing equilibrium points of the system is of interest.
Finally, defining and implementing a stop criterion in case the algorithm converged to stationary points which do not coincide with roots of $\mbf{f}$ would be interesting. 
This study is left by the authors for future work.
    
\begin{acknowledgements}
S.B. was granted a Ph.D. scholarship by Roche S.p.A., Italy. F.B. and S.S. have been partially supported by Gruppo Nazionale per il Calcolo Scientifico (GNCS-INdAM).\\

\noindent
\textbf{Data Availability}
The datasets and codes generated and analysed during the current study are available in the GitHub repository, \url{https://​github.​com/​theMI​DAgro​up/​CRC_​CRN.git}
\end{acknowledgements}

\bibliographystyle{spmpsci_unsrt}
\bibliography{biblio.bib}




                       




\end{document}